\documentclass[11pt]{article}
\pdfoutput=1
\usepackage{array}
\usepackage{cite}
\usepackage{calc}
\usepackage{caption}
\usepackage{graphicx}
\usepackage{epstopdf}
\usepackage{psfrag}
\usepackage{subcaption}
\usepackage{stfloats}
\usepackage{longtable}
\usepackage{amsmath,amsfonts,amssymb,amsbsy,bm,paralist,theorem,ifthen,color}
\usepackage{algorithm}
\usepackage{algpseudocode}
\usepackage[dvipsnames]{xcolor}
\DeclareMathOperator{\argmin}{\mbox{argmin}}

\def\argmin{\mathop{\rm argmin}}

\def\dom{\mbox{dom\,}}

\def\QED{\hfill{\bf Q.E.D.}\smskip}

\newcommand{\beq}{\begin{equation}}
\newcommand{\eeq}{\end{equation}}
\newcommand{\st}{{\rm s.t.}}

\bigskip

\newtheorem{theorem}{Theorem}
\newtheorem{lemma}{Lemma}

\def\smskip{\par\vskip 5 pt}

\def\QED{\hfill{\bf Q.E.D.}\smskip}

\begin{document}
	\title{Gradient-Free Multi-Agent Nonconvex Nonsmooth Optimization}
		\author{{Davood Hajinezhad}
			\thanks{ Department of Mechanical Engineering and Materials Science, Duke University, USA,  Email: \texttt{dhajinezhad@gmail.com}} 
			\and {Michael M. Zavlanos}
			\thanks{Department of Mechanical Engineering and Materials Science, Duke University, USA,  Email: \texttt{michael.zavlanos@duke.edu}}  
		}	
\maketitle
	
\begin{abstract}
	
	In this paper we consider the  problem of minimizing the sum of  nonconvex and possibly nonsmooth functions over a connected multi-agent network, where  the agents have partial knowledge about the global cost function and can only access the zeroth-order information (i.e., the functional values) of their local cost functions. 
	We propose and analyze a distributed primal-dual gradient-free algorithm for this challenging problem. We show that by appropriately choosing the parameters, the proposed algorithm converges to the set of first order stationary solutions with a provable global sublinear convergence rate.
	Numerical experiments demonstrate the effectiveness of our proposed method for optimizing nonconvex and nonsmooth problems over a network.
	
\end{abstract}

\section{Introduction}\label{sub:intro}

Consider a network with $N$ distributed agents that collectively  solve the following optimization problem
\begin{align}\label{sum}
&\min_{x \in \mathbb{R}^M} \; f(x):=\sum_{i=1}^{N} f_i(x).
\end{align}
Here, $f_i: \mathbb{R}^M\to \mathbb{R}$ is  possibly a nonconvex nonsmooth function, that is only available to agent $i$. Such  distributed optimization problems arise in many applications such as machine learning \cite{hong15busmm_spm, davood_NIPS2016_nestt}, resource allocation \cite{bjornson13}, robotic networks \cite{zavlanos2013network}, and signal processing \cite{Cevher14}. See \cite{BoydADMMsurvey2011} for more applications.

Many distributed optimization methods have been proposed to solve problem \eqref{sum}. Many of them rely on consensus between the agents and assume that the cost functions are convex.   One of the first such methods is the distributed subgradient (DSG) algorithm \cite{Nedic09subgradient}. Subsequently, a number of  similar consensus-based algorithms were also proposed to solve distributed convex optimization problems in the form of \eqref{sum}; see, e.g.,   \cite{nedic2009distributed,lobel2011distributed,nedic2015distributed, soomin2016approximate}. These methods only converge to a neighborhood of the solution set unless they use diminishing stepsizes, which often makes them slow. Faster  algorithms using constant stepsizes include the incremental aggregated gradient (IAG) method \cite{gurbuzbalaban2015convergence}, the exact first-order algorithm (EXTRA) \cite{shi2014extra}, and Accelerated Distributed Augmented Lagrangian (ADAL) algorithm. See also  \cite{chatzipanagiotis2015augmented,soomin2018complexity,nicol2016distributed} for optimization of convex problems.

Optimization of nonconvex functions is a much more challenging problem.  Only recently, have there been developed a few nonconvex
distributed optimization algorithms  motivated by applications in resource allocation  in ad-hoc network \cite{tychogiorgos2013non}, sparse PCA \cite{hajinezhad2015nonconvex}, and flow control in communication networks \cite{sun2016distributed}; see also  \cite{bianchi2013convergence,tatarenko2017non,Lorenzo16,sun2016distributed,hong_prox_pda,Chatzipanagiotis2017nonconvex,hajinezhad2018alternating, tarzanagh2017nonmonotone, peyghami2015relaxed} for additional algorithms developed for nonconvex optimization problems.

Regardless of  convexity and/or smoothness,  all the aforementioned methods require that either the first order (gradient/subgradient) information or the explicit form of the objective function is available to the agents.  However, in many important practical situations such information can be expensive to obtain, or even impossible. Examples include,   simulation-based optimization where  the objective function can only be evaluated using repeated simulation \cite{spall2003simulation}, training deep neural networks where 
the relationship between the variables and the cost function is too complicated to derive an explicit form of the gradient \cite{lian2016comprehensive}, and bandit optimization where a player optimizes a sequence of cost functions having only knowledge of a single function value each time  \cite{agarwal2010optimal}. In these cases, the zeroth-order information about the objective function values is often readily available. Such zeroth-order information can be obtained through a stochastic  zeroth-order oracle ($\mathcal{SZO}$). In particular, suppose that $\hat{x}\in\mbox{dom}(f_i)$, then the $i$th $\mathcal{SZO}$ at agent $i$  returns  a noisy  version of $f_i(\hat{x})$ denoted by $\mathcal{H}_i(\hat{x},\xi)$ that satisfies 
\begin{align}\label{eq:bb}
{\mathbb{E}_\xi[\mathcal{H}_i(\hat{x},\xi)] = f_i(\hat{x}),}
\end{align}
where $\xi$ is a random variable representing the noise.

Recently, centralized zeroth-order optimization has received significant attention. In  \cite{nesterov2011random}, Nesterov proposed a general framework for analyzing zeroth-order  algorithms and provided the global convergence rates for both convex and nonconvex problems. 
In \cite{ghadimi2013stochastic}, the authors  established a stochastic zeroth-order method, which again can deal with both convex and nonconvex (but smooth) optimization problems.  In \cite{duchi2015optimal} a Mirror Descent based zeroth-oder algorithm was proposed for solving convex  optimization problems.
Both these  zeroth-order algorithms are centralized and cannot be implemented over a multi-agent network. A few recent works  \cite{yuan2015randomized,yuan2015gradient,yuan2016zeroth} considered zeroth-order distributed convex (possibly nonsmooth) problems, but  none of these works can address nonconvex problems. For distributed nonconvex but smooth optimization problem using zeroth-order information a primal-dual algorithm has been proposed in \cite{hajinezhad2017zeroth, davood2019zone}. 

In this paper, we propose a new algorithm for distributed nonconvex and nonsmooth optimization with zeroth-order information. Specifically, we first show that this problem can be reformulated as a linearly constrained optimization problem over a connected multi-agent network. Then,  we propose a nonconvex primal-dual based algorithm, which requires only  local communication among the agents, and utilizes local zeroth-order information.  
Theoretically, We show that the current solution converges approximately to a stationary solution of the problem. We also provide numerical results that corroborate the theoretical findings. 

\noindent {\bf Notation.} 
Given a vector $a$ and a matrix $A$, we use $\|a\|$ and $\|A\|$ to denote the Euclidean norm of vector $a$, and  spectral norm of matrix $A$, respectively. $A^\top$ represents the transpose of matrix $A$. We define $\|a\|^2_A:=a^TAa$. The notation {$\langle a,b\rangle$} is used to denote the inner product of two vectors $a$, $b$. For matrices $A$ and $B$, $A\otimes B$ is the Kronecker product of $A$ and $B$.   To denote an $M\times M$ identity matrix we use $I_M$. $\mathbb{E}[\cdot]$ denotes the expectation with respect to all random variables, and $\mathbb{E}_v[\cdot]$ denotes the expectation with respect to the random variable $v$.

\section{ Problem Definition and Proposed Algorithm} 
Let us define a graph $\mathcal{G}= \{\mathcal{V}, \mathcal{E}\}$, where $\mathcal{V}$ is the node set with $|\mathcal{V}| = N$, and $\mathcal{E}\subseteq \mathcal{V}\times\mathcal{V}$ is the edge set  with $|\mathcal{E}|=E$. We assume that $\mathcal{G}$ is undirected, meaning that if $(i,j)\in\mathcal{E}$ then $(j,i)\in\mathcal{E}$. Moreover, every agent $i$ can only communicate with its direct neighbors in the set $\mathcal{N}_i=\{j\in\mathcal{V}\; |\; (i,j)\in\mathcal{E}\}$, and let $d_i=|\mathcal{N}_i|$ denote the degree of node $i$. We assume that the graph $\mathcal{G}$ is connected, meaning that there is a path, i.e., a sequence of nodes where consecutive nodes are neighbors, between any two nodes in $\mathcal{G}$.

In order to decompose problem \eqref{sum} let us introduce $N$ new variables $x_i\in\mathbb{R}^M$ that are local to every agent $i$. Then, problem \eqref{sum} can be reformulated as follows: 
\begin{align}\label{pr-due}
\min_{\{x_i\}}&\quad \sum_{i=1}^{N}f_i(x_i),\\
\st &\quad x_i = x_j\; \forall~(i,j)\in\mathcal{E}. \nonumber
\end{align}
The set of constraints $x_{i}=x_j$ enforce consensus on the local variables $x_i$ and $x_j$ for all neighbors $j\in\mathcal{N}_i$. We stack all the local variables $x_i$ in a vector 
$x: = \{x_i\}\in\mathbb{R}^{Q\times 1}$, where   $Q=NM$.  
Moreover, we define the Degree matrix $\tilde{D}\in\mathbb{R}^{N\times N}$  to be a diagonal matrix where $\tilde{D}(i,i)=d_i$; let $D =  \tilde{D}\otimes I_M\in\mathbb{R}^{Q\times Q}$. 
For a given graph $\mathcal{G}$, the { incidence matrix} $\tilde{A}\in \mathbb{R}^{E\times N}$ is a matrix where  $\tilde{A}(k,i)=1$ and $\tilde{A}(k,j)=-1$, where $(i,j)$ is the $k$th edge of $\mathcal{G}$; the rest of the entries of $\tilde{A}$ are all zero. Let {$A=\tilde{A}\otimes I_M\in\mathbb{R}^{EM\times Q}$}.
Finally, we define the Signed and Signless Laplacian matrix, denoted by $L^-$ and $L^+$, respectively as
\begin{align}
L^{-}&:=A^\top A\in\mathbb{R}^{Q\times Q},\label{eq:slap}\\
L^{+}&:=2D-A^\top A\in\mathbb{R}^{Q\times Q}.\label{eq:sllap}
\end{align}  
Using the above notations,  problem \eqref{pr-due} can be written  in the following compact form:
\begin{align}\label{pro:compact}
\min_{x\in \mathbb{R}^Q} \; f(x),\quad \st \; \;  Ax = 0. 
\end{align}

\subsection{ Preliminaries}  
In this section, we first introduce some standard techniques presented in \cite{nesterov2011random} for approximating and smoothing the gradient of a given function. Suppose that $\phi\in \mathbb{R}^Q$ is a Gaussian random vector and let $\mu>0$ be some smoothing parameter. The smoothed version of function $f$ is defined as
\begin{align}\label{eq:smooth}
f_{\mu}(x) &= \mathbb{E}_\phi(f(x+\mu \phi))=\frac{1}{(2 \pi)^{\frac{Q}{2}}} \int f(x+\mu \phi) e^{-\frac{1}{2}\|\phi\|^2} \,d\phi.
\end{align}
Then it can be shown that the function $f_{\mu}$ is  differentiable  and its gradient is given by Eq. (22) in \cite{nesterov2011random}
\begin{align} \label{eq:grad_smooth_ver}
\hspace{-.3cm}\nabla f_{\mu}(x) = \frac{1}{(2 \pi)^{\frac{Q}{2}}} \int \frac{f(x+\mu \phi)-f(x)}{\mu} \phi e^{-\frac{1}{2}\|\phi\|^2} \,d\phi.
\end{align}
Further, assuming that the original function $f$ is Lipschitz continuous, i.e., there exists $L_0$ such that $|f(x)-f(y)|\leq L_0|x-y|$ for all $x,y\in\dom(f)$,  it can be shown that (see \cite[Lemma 2]{nesterov2011random}) $\nabla f_\mu$ is also Lipschitz continuous with constant $L_1=\frac{2L_0\sqrt{Q}}{\mu}$. In other words, for all $x,y\in\dom(f_\mu)$ we have
\begin{align}\label{eq:smooth_lip}
\|\nabla f_\mu(x)-\nabla f_\mu(y)\|\leq L_1\|x-y\|.
\end{align}

Let $\mathcal{H}(x,\xi)$ denote the noisy functional value of the  function $f$ obtained from an associated $\mathcal{SZO}$ as in equation \eqref{eq:bb}. In view of  \eqref{eq:grad_smooth_ver}, the gradient of $f_\mu(x)$   can be approximated as
\begin{align}\label{eq:grad:estimate}
G_\mu(x, \phi, \xi) = \frac{\mathcal{H}(x+\mu\phi, \xi) - \mathcal{H}(x, \xi)}{\mu}\phi,
\end{align}
where the constant  $\mu>0$ is the smoothing parameter.  
It can be easily checked that  $G_\mu(x, \phi, \xi)$ is an unbiased estimator of $\nabla f_\mu (x)$, i.e.,
\begin{align} \label{eq:def_ex_G_mu}
\hspace{-.33cm}\mathbb{E}_{\xi, \phi}[G_{\mu}(x, \phi, \xi)] & = \mathbb{E}_\phi \left[ \mathbb{E}_{\xi}[G_{\mu}(x, \phi, \xi)\mid \phi] \right]= \nabla f_{\mu}(x).
\end{align}
For simplicity we define $\zeta:=(\xi,\phi)$. For a given number $J$ of independent samples of $\{\zeta_j\}_{j=1}^J$, we define  the sample average $\bar{G}_{\mu}(x,\zeta):=\frac{1}{J}\sum_{j=1}^J G_\mu(x, \zeta_j)$, where $\zeta:=\{\zeta_j\}_{j=1}^J$. 
It is easy to see that for any $J\ge 1$,  $\bar{G}_{\mu}(x,\zeta)$ is also an unbiased estimator of $\nabla f_\mu (x)$.

\subsection{ The Proposed Algorithm}
In this part we propose a primal-dual algorithm for the distributed optimization problem \eqref{pro:compact}. 
Let $\lambda_{ij}\in\mathbb{R}^M$ be the multiplier associated with the consensus constraint $x_i -x_j= 0$ for each $(i,j)\in\mathcal{E}$. Moreover,  stack all $\lambda_{ij}$'s in a vector $\lambda=\{\lambda_{ij}\}_{(i,j)\in\mathcal{E}} \in\mathbb{R}^{EM}$.  Then, the augmented Lagrangian (AL) function for problem \eqref{pro:compact} is given by
\begin{equation}\label{eq:lag}
U_\rho(x,\lambda) := f(x)+\langle\lambda,Ax\rangle+ \frac{\rho}{2}\left\|Ax\right\|_2^2.
\end{equation}
where $\rho>0$ is a constant. Moreover, as in  \eqref{eq:smooth} define the smoothed version $f_{i,\mu}$ of the local function $f_i$.  At iteration $r$ of the algorithm we obtain an unbiased estimation of the gradient of local function $f_{i,\mu}(x_i^r)$ as follows. For every sample $j\in\{1,2,\cdots,J\}$ we generate a random vector $\phi^r_{i,j}\in\mathbb{R}^M$  from an i.i.d standard Gaussian distribution and calculate $\bar{G}_{\mu,i}(x_i^r,\zeta^r_{i})\in\mathbb{R}^{M}$ similar to \eqref{eq:grad:estimate} by
\vspace{-.2cm}
\begin{align}\label{eq:grad:approx}
&{ \bar{G}_{\mu,i}(x_i^r,\zeta^r_{i}})=\frac{1}{J}\sum_{j=1}^J\frac{\mathcal{H}_i(x_i^r+\mu\phi^r_{i,j}, \xi^r_{i,j}) - \mathcal{H}_i(x_i^r, \xi^r_{i,j})}{\mu}\phi^r_{i,j}.
\end{align}
Define ${  G^{J,r}_{\mu}}:=\{\bar{G}_{\mu,i}(z_i^r,\zeta^r_{i})\}_{i=1}^N\in\mathbb{R}^{Q}$. The following theorem bounds the norm of $G^{J,r}_{\mu}$. 
\begin{theorem}[Theorem 4 \cite{nesterov2011random}]
	If $f$ is a Lipschitz continuous function with constant $L_0$, then
	\begin{align}\label{eq:nes_them4}
	\mathbb{E}_{\zeta}[\|{G}_{\mu}^{J,r}\|^2]\leq \frac{L_0^2(Q+4)^2}{J^2}.
	\end{align}
\end{theorem}
\begin{algorithm}[tb]
	\caption{The Proposed Algorithm for problem \eqref{sum}}
	\label{alg}
	\begin{algorithmic}[1]
		\State {\bfseries Input:} Degree matrix $D\in\mathbb{R}^{Q\times Q}$, total number of iterations $T\ge 1$, number of samples $ J\ge1$,  smoothing parameter $\mu>0$
		\State  {\bfseries Initialize:} Primal variable ${x}^0\in\mathbb{R}^{Q}$, dual variable {$\lambda^0=0\in\mathbb{R}^{EM}$} 
		\For{$r=0$ {\bfseries to} ${T-1}$} 
		\Statex Update primal variable $x$ and dual variable $\lambda$ by 
		\begin{align}
		x^{r+1}& = \argmin_{x} ~\langle {G^{J,r}_{\mu}}+A^\top\lambda^r+\rho A^\top Ax^r,x-x^r\rangle +\rho\|x-x^r\|^2_{D},\texttt{(Primal Step)}\label{eq:x:zero}\\
		\lambda^{r+1}& = \lambda^r +\rho A x^{r+1}. \qquad\texttt{ (Dual Step)}\label{eq:lambda:zero}
		\end{align}
		\EndFor
		\State { Choose uniformly randomly $u\in\{0,1,\cdots,T-1\}$}
		\State {\bfseries Output:} $(z^u, \lambda^u)$.
	\end{algorithmic}
\end{algorithm} 
Our proposed algorithm is summarized in Algorithm \ref{alg}. In the primal step \eqref{eq:x:zero}, an {\it approximate gradient descent} step is taken towards minimizing the augmented Lagrangian function with respect to $x$. 
In particular, in the first-order approximation of the AL function the true gradient of the function $f(x^r)$ is approximated by the noisy zeroth-order estimate $ G_{\mu}^{J,r}$ and then, a matrix-weighted quadratic penalty $\rho\|x-x^r\|_D$ is used. This term is critical for the algorithm itself, as well as for the analysis. The dual step \eqref{eq:lambda:zero} is then performed, which is a {\it gradient ascent} step over the dual variable $\lambda$.

To see how Algorithm \ref{alg} can be implemented in a distributed way, consider the optimality condition for \eqref{eq:x:zero} as
\begin{align}\label{eq:kkt:x:compact}
G_\mu^{J,r}+A^\top(\lambda^r+\rho Ax^r)+2\rho D(x^{r+1}-x^r)=0.
\end{align}

Utilizing \eqref{eq:slap} and \eqref{eq:sllap}, we have
\begin{align*}
x^{r+1}=\frac{1}{2\rho}D^{-1}\bigg[\rho L^{+}x^r-G_\mu^{J,r}+A^\top\lambda^r\bigg].
\end{align*}

To implement this primal iteration, each agent $i$ only requires local information as well as information from its neighbors $\mathcal{N}_i$. This is because $D$ is a diagonal matrix and the structure of the matrix $L^+$  ensures that the $i$th block vector of  $L^{+} x^r$ is only related to $x^{r}_j, \; j\in\mathcal{N}_i$. For the dual step w.l.o.g we assign the dual variable $\lambda_{ij}$  to node $i$ and therefore, from \eqref{eq:lambda:zero} we have
\begin{align}
\lambda_{ij}^{r+1}=\lambda_{ij}^r + \rho(x_i^{r+1}-x_j^{r+1}),
\end{align} 
which only requires the local information as well as  information from the neighbors in $\mathcal{N}_i$.   

\section{ The Convergence Analysis}\label{sec:analysis}
In this section we study the convergence  of  Algorithm \ref{alg}.  We make the following assumptions. 

\noindent{\bf Assumptions  A.} We assume that
\begin{itemize}
	\item[A1.]  The function $f$ is Lipschitz continuous.
	\item[A2.]  The function $f$ is lower bounded. 
\end{itemize}

The above assumptions on the objective $f$ are quite standard  in the analysis of first order optimization Algorithms (1).  
To simplify notation let $\mathcal{F}^r:=\sigma(\zeta_1,\zeta_2,\cdots\zeta_r)$ be the $\sigma$-field generated by the entire history of algorithm up to iteration $r$, $\sigma_{\min}$ be the smallest nonzero eigenvalue of $A^\top A$, and $w^r:=(x^{r+1}-x^r) -(x^{r}-x^{r-1})$ be the successive difference of the differences of the primal iterates. 
In the analysis that follows we will make use of the following relations:
\begin{itemize}
	\item For any given vectors $a$ and $b$ we have
	\vspace{-.2cm}
	{
		\begin{align}
		&\langle b-a, b\rangle = \frac{1}{2}(\|b\|^2+\|a-b\|^2-\|a\|^2), \label{eq:rel1}\\
		& \langle a,b\rangle \le \frac{1}{2\epsilon}\|a\|^2 + \frac{\epsilon}{2}\|b\|^2; \quad \forall~ \epsilon>0. \label{eq:rel3}
		\end{align}}
	\item For $n$ given vectors  $a_i$ we have that
	\vspace{-.2cm}
	{
		\begin{align}\label{eq:rel2}
		\bigg\|\sum_{i=1}^na_i\bigg\|^2\le n\sum_{i=1}^{n}\big\|a_i\big\|^2.
		\end{align}}
\end{itemize}

Our proof consists of a series of lemmas leading to the main convergence rate result.
In our presentation we try to provide insights about the analysis steps, while the proofs of the results are relegated to the appendix. 

First we bound the difference between $\nabla f_{\mu}(x^r)$ and its unbiased estimation $G^{J,r}_{\mu}$ as follows:
\begin{align}\label{eq:bd_unbised}
\mathbb{E}_\zeta\|G^{J,r}_{\mu}-\nabla f_\mu(x^r)\|^2 &= \mathbb{E}_\zeta\|G^{J,r}_{\mu}-\mathbb{E}_\zeta[G^{J,r}_{\mu}]\|^2\nonumber\\
&\leq \mathbb{E}_\zeta\|G^{J,r}_{\mu}\|^2\nonumber\\
&\leq \frac{L_0^2(Q+4)^2}{J^2},
\end{align}
where the last inequality follows from \eqref{eq:nes_them4}. The next lemma bounds the change of the dual variables  by that of the primal variables. The proofs of the results that follow can be found in the appendix.
\begin{lemma}\label{lemma:mu:bound}
	Suppose Assumptions A hold true. Let $L_1$ denote the gradient Lipschitz constant for function $f_\mu$. Then we have the following inequity:
	{
		\begin{align}\label{eq:mu:difference:bound}
		\frac{1}{\rho}\mathbb{E}\|\lambda^{r+1}-\lambda^{r}\|^2&\leq \frac{9L_0^2(Q+4)^2}{\rho\sigma_{\min}J^2}+\frac{6L^2_1}{{\rho\sigma_{\min}}}\mathbb{E}	\|x^r-x^{r-1}\|^2+\frac{3\rho\|L^{+}\|}{\sigma_{\min}}\mathbb{E}\|w^r\|_{L^{+}}^2.
		\end{align}}
\end{lemma}

The next step is the key in our analysis. We define the smoothed version of the AL function in a similar way as \eqref{eq:smooth} and denote it by $U_{\rho,\mu}(x,\lambda)$. For notational simplicity let us define $U^{r+1}_{\rho,\mu}:=U_{\rho,\mu}(x^{r+1},\lambda^{r+1})$. From equation \eqref{eq:smooth_lip} we know that function $f_{\mu}$ is Lipschitz continuous with constant $L_1$. Now let $c>0$ be some positive constant and set    $k:=2\big(\frac{6{L_1}^2}{\rho\sigma_{\min}}+\frac{3c{L_1}}{2}\big)$.  Moreover, we define $$V^{r+1}:=\frac{\rho}{2}\big(\|Ax^{r+1}\|^2+\|x^{r+1}-x^r\|^2_{B}\big),$$ where $B:=L^++\frac{k}{c\rho}I_Q$. Finally we define the following potential function:
\begin{align}\label{eq:pot_def}
P^{r+1}:=U^{r+1}_{\rho,\mu}+ cV^{r+1}.
\end{align}
We study the behavior of the proposed potential function as the algorithm proceeds.   
\begin{lemma}\label{lem:bd:pot}
	Suppose Assumptions A hold true. We have that
	\begin{align}
	\mathbb{E}\big[P^{r+1}-P^{r}\big]&\le -\alpha_1\mathbb{E}\|x^{r+1} - x^r\|^2 -\alpha_2\mathbb{E}\|w^r\|_{L^{+}}^2 + \alpha_3\frac{L_0^2(Q+4)^2}{J^2},\label{eq:diff:bd}
	\end{align}
	where
	{  
		\begin{align}\label{eq:c}
		&\alpha_1:= \rho^2-(2cL_1+L_1/2+L_1^2/2+1/2)\rho-\frac{6L_1^2}{\sigma_{\min}},\nonumber\\
		& \alpha_2:=\frac{3\rho\|L^{+}\|}{\sigma_{\min}}-\frac{c \rho}{2},\; \alpha_3=\frac{9}{\rho\sigma_{\min}}+\frac{6c+1}{L_1}.
		\end{align}}
\end{lemma}

Note that the constants $\alpha_1$ and $\alpha_2$ in \eqref{eq:diff:bd} can be made positive as long as we choose constants $c$ and $\rho$ large enough. In particular, the following conditions are sufficient to have positive $\alpha_1,\alpha_2$
\begin{align}
&c>\frac{6\|L^{+}\|}{\sigma_{\min}},\quad\rho>b+\sqrt{b^2+6L_1^2/\sigma_{\min}},\label{eq:bd:c}
\end{align}where $b=(cL_1+L_1/4+L_1^2/4+1/4)$. 

The key insight obtained from Lemma \ref{lem:bd:pot} is that, a proper combination of a primal objective (i.e., the AL function), and the dual gap (i.e., the violation of the feasibility) can be served as the potential function that guides the progress of the algorithm. 

In the next lemma we show that $P^{r+1}$ is lower bounded.
\begin{lemma}\label{lemma:lower:bound}
	Suppose Assumptions A hold true, and constant $c$ is selected as  $c\geq \frac{2\|L^+\|}{\sigma_{\min}}.$
	Then there exists a constant $\underbar{P}$ that is independent of the total number of iterations $T$ so that
		\begin{align}\label{eq:lower:bound}
		\mathbb{E}[P^{r+1}]\ge \underline{P}>-\infty, \quad \forall~r\geq 1.
		\end{align}
\end{lemma}

To characterize the convergence rate of Algorithm \ref{alg}, let us define the {\it stationarity gap} of the smoothed version of problem \eqref{pro:compact} as
\begin{align}\label{eq:opt_gap}
\hspace{-.3cm}\Phi_\mu(x^r,\lambda^{r-1}): = \mathbb{E}\left[\|\nabla_{x}U_{\rho,\mu}(x^r, \lambda^{r-1})\|^2+\|A x^r\|^2\right].
\end{align}
It can be easily checked that $\Phi_\mu(x^*,\lambda^*)= 0$ if and only if $(x^*, \lambda^*)$ is a KKT point of the smoothed version of problem \eqref{pro:compact}. For simplicity let us denote $\Phi_\mu^{r}:=\Phi_\mu(x^r, \lambda^{r-1})$.

At this point we are ready to combine the previous results to obtain our main theorem.
\begin{theorem}\label{thm:conv}
	Suppose Assumptions A hold true,   the penalty parameter $\rho$ satisfies  the condition given in Lemma \ref{lem:bd:pot}, and the constant $c$ satisfies  {$c\geq \frac{6\|L^{+}\|}{\sigma_{\min}}.$}
	Then, there exist constants $\gamma_1, \gamma_2>0$ such that
		\begin{align}
		\mathbb{E}_u[\Phi_\mu^u]\le \frac{\gamma_1}{T} + \gamma_2\frac{L_0(Q+4)^2}{J^2}.
		\end{align}
\end{theorem}

From Theorem \ref{thm:conv} we can observe that there exists always a constant term in the right-hand-side of the stationarity gap. Therefore, no matter how many iterations we run the algorithm, we always converge to a neighborhood of a stationary  point. However, if we choose the number of samples $J\in\mathcal{O}({\sqrt{T}})$, we have the following bound: 
\begin{align}
\mathbb{E}_u[\Phi^u_\mu]\le \frac{\gamma_1}{T} +\frac{\gamma_2L_0(Q+4)^2}{T},
\end{align}
which verifies the sublinear convergence rate for the algorithm.

\section{ Numerical Results}
In this section we illustrate the proposed algorithm through numerical simulations. For our experiments we study two nonconvex distributed optimization problems. 

\begin{figure}[t!]
	\centering
	\begin{subfigure}{.49\textwidth}
	\includegraphics[width=2.7 in]{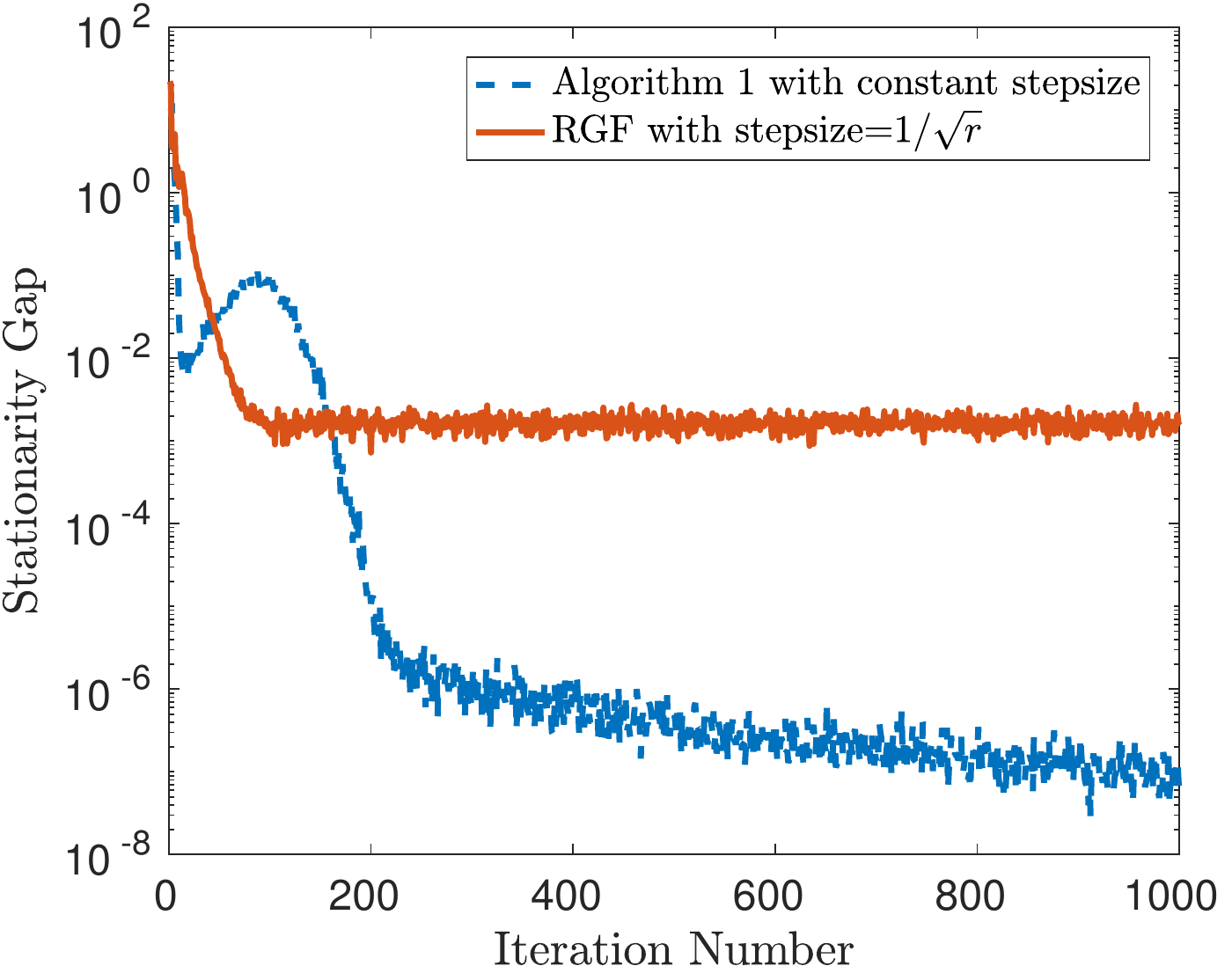}
	\vspace{-0.2cm}
	\caption{\footnotesize  Comparison  of proposed Algorithm 1, and RGF algorithm \cite{yuan2015randomized} in terms of the stationarity gap for the nonconvex nonsmooth distributed optimization problem \eqref{p1_problem}. } 
	\label{fig:opt_gap:cons}
	\end{subfigure}
\hfill
\begin{subfigure}{.49\textwidth}
	\centering
	\includegraphics[width=2.7 in]{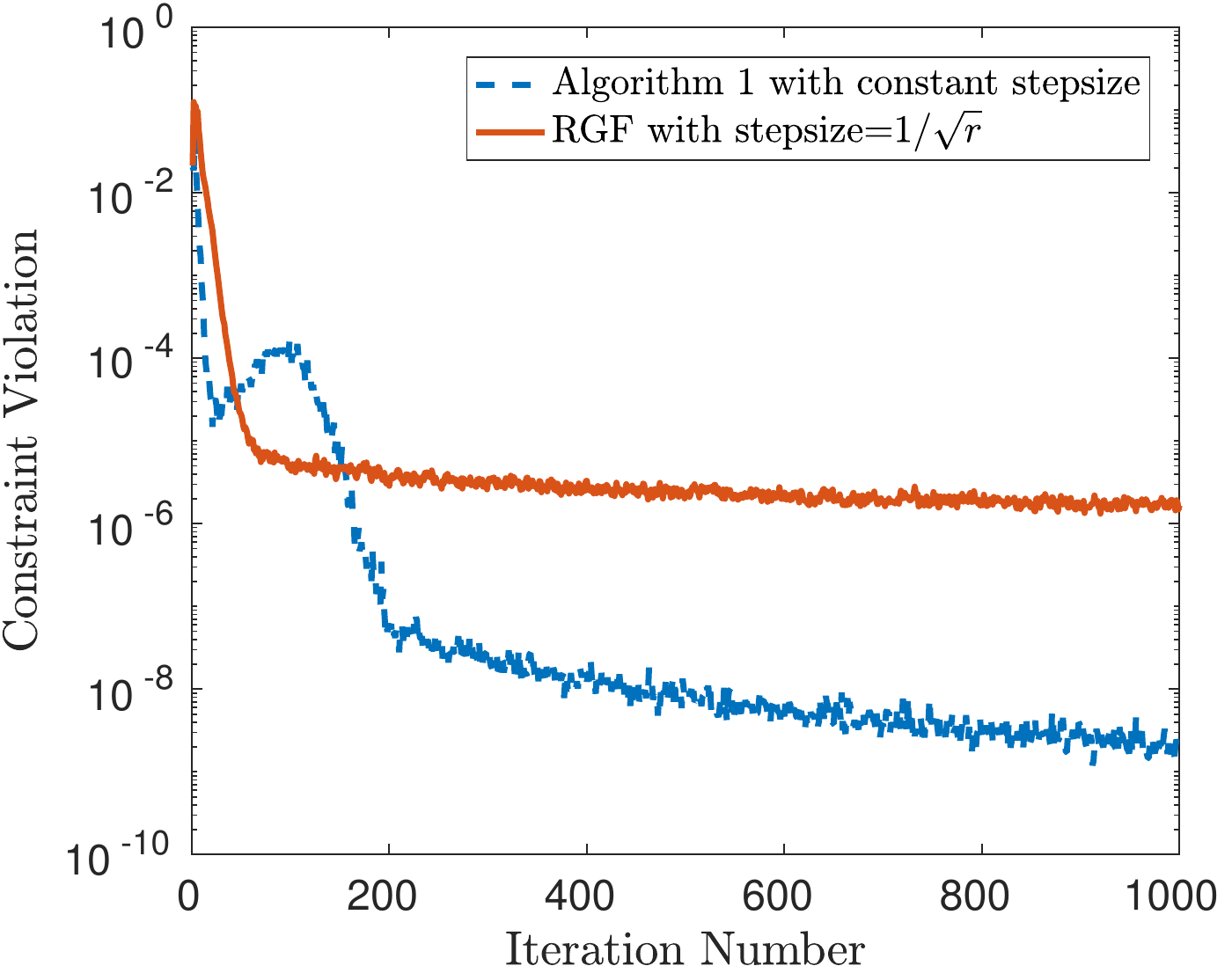}
	\vspace{-0.2cm}
	\caption{\footnotesize Comparison  of proposed Algorithm 1, and RGF algorithm \cite{yuan2015randomized} in terms of the constraint violation (i.e. $\|Ax\|$) for the nonconvex nonsmooth distributed optimization problem \eqref{p1_problem}. }
	\label{fig:cons_vio:cons}
\end{subfigure}   
\end{figure}

\begin{figure}[t!]
	\centering
	\begin{subfigure}{.49\textwidth}
	\includegraphics[width=2.6 in]{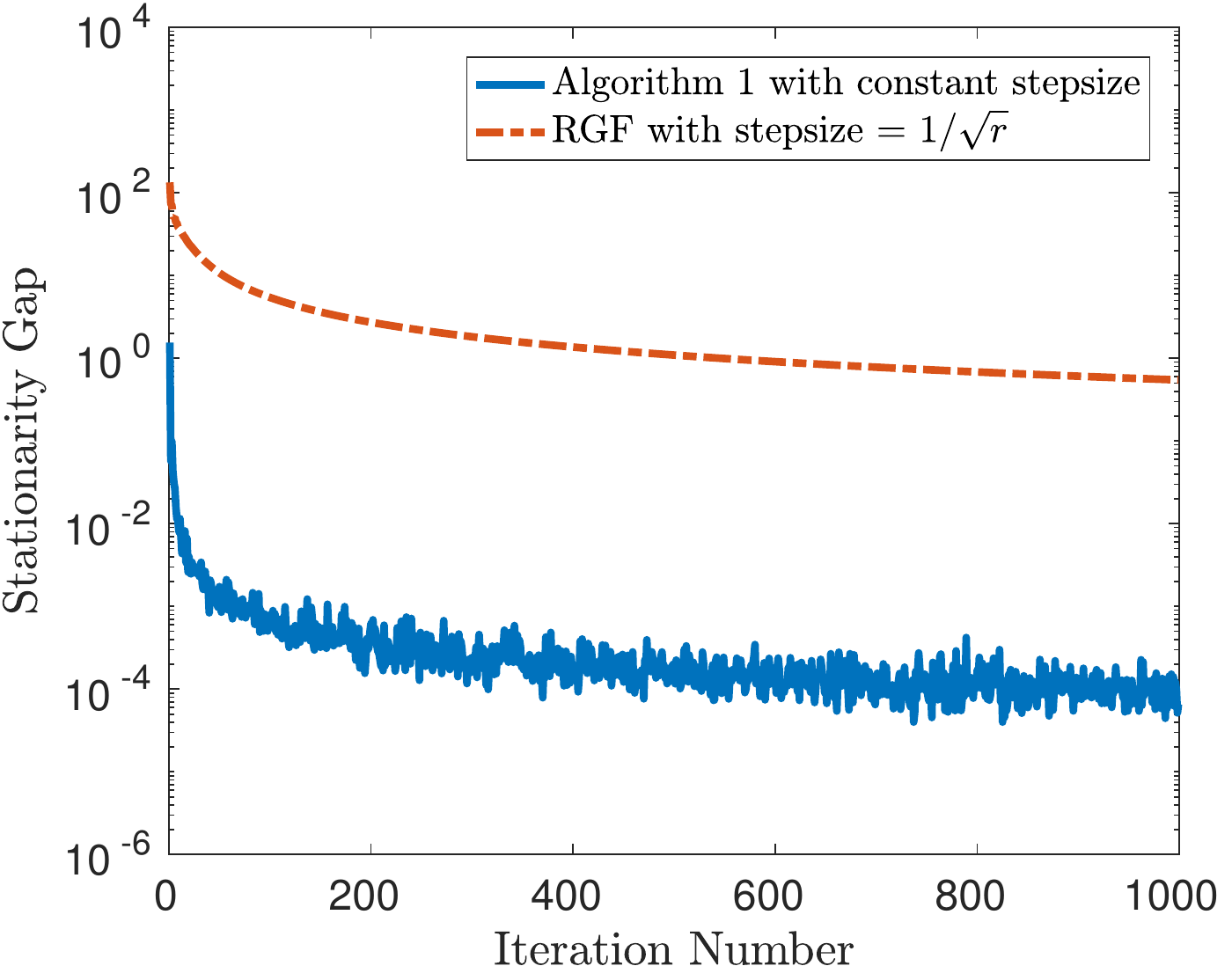}
	\caption{\footnotesize  Comparison  of proposed Algorithm 1, and RGF algorithm \cite{yuan2015randomized} in terms of the stationarity gap for binary classification problem.}
	\label{fig:opt_gap_bin}
\end{subfigure}
\begin{subfigure}{.49\textwidth}
	\centering
	\includegraphics[width=2.6 in]{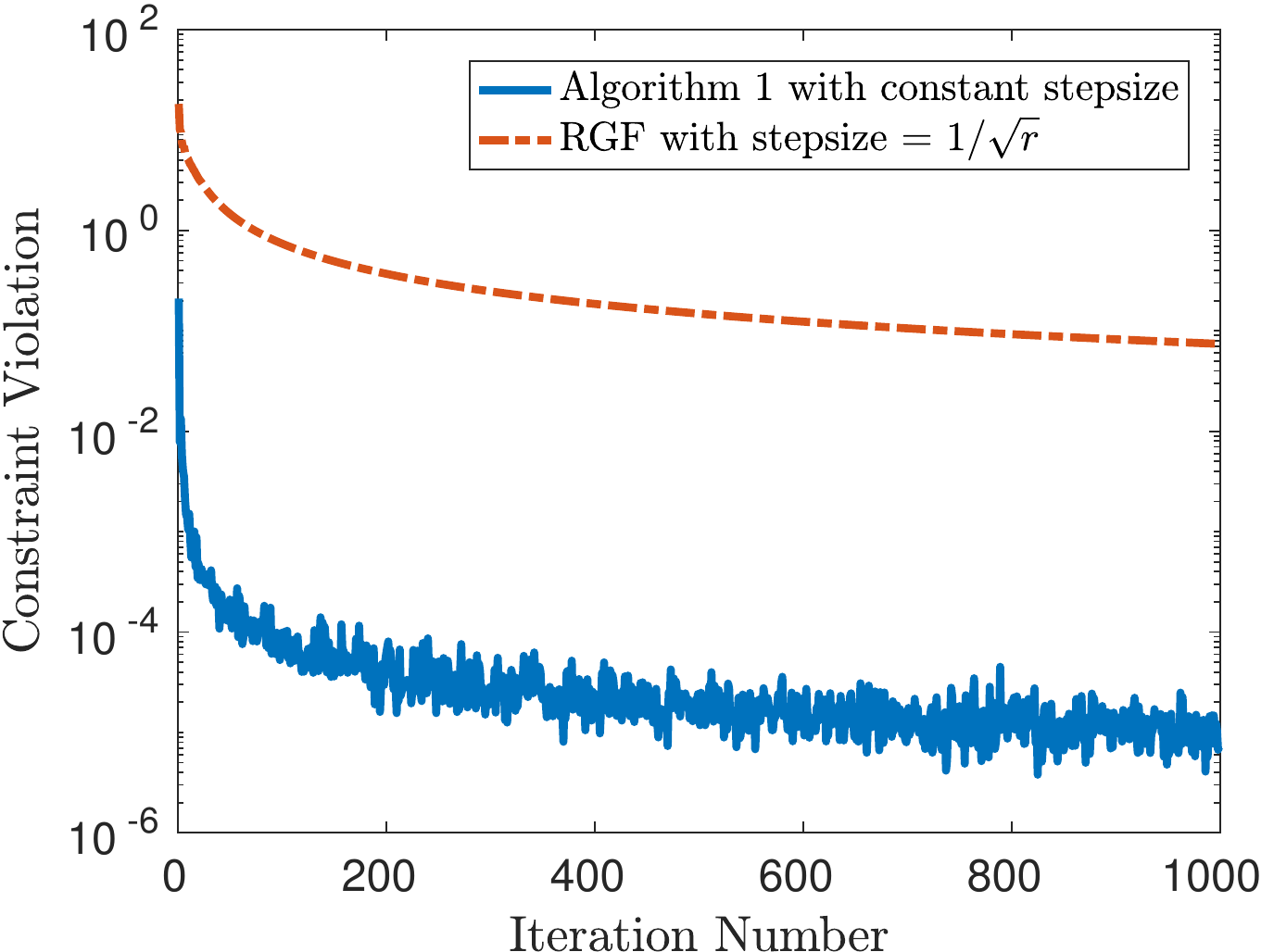}
	\caption{\footnotesize  Comparison  of proposed Algorithm 1, and RGF algorithm \cite{yuan2015randomized} in terms of the constraint violation (i.e. $\|Ax\|$) for binary classification problem.}
	\label{fig:cons_vio_bin}
\end{subfigure}
\end{figure}

First, we consider a simple nonconvex nonsmooth distributed optimization problem defined bellow:
\begin{equation} \label{p1_problem}
\mathop {\min }\limits_{x\in \mathbb{R}^Q} \sum\limits_{i = 1}^N {{f_i}({x_i})},\quad \textrm{s.t.}\;x_i=x_j,\; \forall (i,j)\in\mathcal{E},
\end{equation} 
where for each agent we have $${f_i}({x_i}) = |\cos(x_i)+\|x_i\|+\exp(x_i)|.$$ In \eqref{p1_problem} the problem dimension $M=1$ and the number of nodes in the network is $N=10$. Therefore, $Q=MN=10$.  The details of the underlying graph are discussed in \cite{YildizScag08}. We compare Algorithm \ref{alg} using a constant stepsize that satisfies \eqref{eq:bd:c} and the Randomize Gradient Free (RGF) algorithm proposed in \cite{yuan2015randomized} using diminishing stepsize $\frac{1}{\sqrt{r}}$. Note that in theory RGF only works for the convex problems. However, we include it here for  the purpose of comparison only. We compare the two algorithms in terms of the stationarity gap defined in \eqref{eq:opt_gap} and the constraint violation $\|Ax\|$. The stopping criterion is set to $T=1000$ iterations and the results are the average over $30$ independent trials. Figs \ref{fig:opt_gap:cons} and \ref{fig:cons_vio:cons} show our comparative results. We observe that the stationarity gap and the consensus error vanish faster for our proposed algorithm than for the RGF algorithm. 
Next, we study a mini-batch binary classification problem using nonconvex nonsmooth regularizers, where each node stores $b$ (batch size) data points. For this problem the local function is given by
	\begin{align*}
	&f_i(x_i) = \frac{1}{Nb}\bigg[\sum_{j=1}^{b}\log (1 + \textrm{exp}( - {y_{ij}}x_i^T{v_{ij}}))+\alpha\log(\epsilon+\|x_i\|_1)\bigg],
	\end{align*}
where $v_{ij} \in \mathbb{R}^M$ and $y_{ij}\in\{1,-1\}$ are the feature vector and the label for the $j$th data point of $ith$ agent \cite{Antoniadis2009}. The nonconvex nonsmooth regularization term $\log(\epsilon+\|x_i\|_1)$ imposes sparsity to vector $x_i$, the constant $\alpha$ controls the sparsity level, and $\epsilon>0$ is a small number.  The network has $N=15$ nodes and each node contains randomly generated  $b=100$ data point. Algorithm \ref{alg} and RGF run for $T=1000$ iterations and Figures \ref{fig:opt_gap_bin} and \ref{fig:cons_vio_bin} illustrate the stationarity gap and the constraint violation versus the iteration counter.  From these plots we can observe that, as in the previous problem, Algorithm \ref{alg} is faster than the RGF. Note again that, in theory, RGF is designed for convex problems only. To the best of our knowledge, our algorithm is the first provable distributed zeroth-order method for nonconvex and nonsmooth problems.

\section{Conclusion}
In this work, we proposed  a distributed gradient-free optimization algorithm to solve nonconvex and nonsmooth problems  utilizing local zeroth-order information. We rigorously analyzed the convergence rate of the proposed algorithm and demonstrated its performance via simulation. To the best of our knowledge, this is the first distributed framework for the solution of nonconvex and nonsmooth distributed  optimization problems that also has a provable sublinear convergence rate. The proposed framework can be used to solve a variety of problems where access to first or second order information is very expensive or even impossible.



\section*{APPENDIX}
 \noindent\subsection{ Proof of Lemma \ref{lemma:mu:bound}}
	From equation  \eqref{eq:lambda:zero} we have
		\begin{align}\label{eq:due:rearg}
		\lambda^{r+1}-\lambda^r=\rho Ax^{r+1}.
		\end{align}
		Equation  \eqref{eq:due:rearg} implies that $\lambda^{r+1}-\lambda^r$ lies in the column space of $A$, therefore we have
		\begin{align}\label{eq:bd:sig}
		\sqrt{\sigma_{\min}}\|\lambda^{r+1}-\lambda^r\|\leq\|A^\top (\lambda^{r+1}-\lambda^r)\|,
		\end{align}
		where $\sigma_{\min}$ denotes the smallest non-zero eigenvalue of $A^\top A$. 
	Utilizing  equation \eqref{eq:due:rearg} and equation \eqref{eq:kkt:x:compact}, we obtain
		\begin{align}\label{eq:kkt:z:app}
		G_{\mu}^{J,r}+A^\top\lambda^{r+1}+\rho L^{+}(x^{r+1}-x^r)=0.
		\end{align}
		Replacing $r$ with $r-1$ in equation \eqref{eq:kkt:z:app} and then using the definition of $w^r:=(x^{r+1}-x^r) -(x^{r}-x^{r-1})$ we obtain
		\begin{align}\label{eq:lam}
		\frac{1}{\rho}\big\|\lambda^{r+1}-\lambda^{r}\big\|^2&\le \frac{1}{{\rho\sigma_{\min}}}\big\|G_{\mu}^{J,r}-G_{\mu}^{J,r-1}{ +}\rho L^{+}w^r\big\|^2\nonumber\\
		&=\frac{1}{\rho\sigma_{\min}}\|G_{\mu}^{J,r}-G_{\mu}^{J,r-1}+\nabla f_{\mu}(x^r)-\nabla f_{\mu}(x^r){ +}\rho L^{+}w^r\|^2\nonumber\\
		&\le \frac{3}{\rho\sigma_{\min}}\|G_{\mu}^{J,r}-\nabla f_{\mu}(x^r)\|^2+\frac{3}{{\rho\sigma_{\min}}}\|\nabla f_{\mu}(x^r)- G_{\mu}^{J,r-1}\|^2+\frac{3\rho^2}{\rho\sigma_{\min}}\| L^{+}w^r\|^2,
		\end{align}
		where the last inequality follows from \eqref{eq:rel2}. Adding and subtracting  $\nabla f_{\mu}(x^{r-1})$ to the second term on the r.h.s of \eqref{eq:lam} and taking the  expectation on both sides gives
	
		\begin{align}\label{eq:bd:dual}
		\frac{1}{\rho}\mathbb{E}\|\lambda^{r+1}-\lambda^{r}\|^2&\le \frac{3}{\rho\sigma_{\min}}\mathbb{E}\|G^{J,r}_{\mu}-\nabla f_{\mu}(x^r)\|^2 \nonumber\\
		&+\frac{6}{\rho\sigma_{\min}}\mathbb{E}\|\nabla f_\mu(x^r)-\nabla f_\mu(x^{r-1})\|^2+\frac{3\rho^2}{\rho\sigma_{\min}}\mathbb{E}\| L^{+}w^r\|^2
		\nonumber\\
		&+\frac{6}{\rho\sigma_{\min}}\mathbb{E}\|\nabla f_{\mu}(x^{r-1})- G_{\mu}^{J,r-1}\|^2\nonumber\\
		&\le \frac{9L_0^2(Q+4)^2}{\rho\sigma_{\min}J^2}+\frac{6L^2_1}{\rho\sigma_{\min}}\mathbb{E}	\|x^r-x^{r-1}\|^2+\frac{3\rho\|L^{+}\|}{\sigma_{\min}}\mathbb{E}\|w^r\|_{L^+}^2,	
		\end{align}
		where the last inequity is true because of \eqref{eq:bd_unbised}, the fact that $\nabla f_\mu(z)$ is  gradient Lipschitz with constant $L_1$,  and the inequality $\|L^+w^r\|^2\le \|L^+\|\|w^r\|^2_{L^{+}}$. The proof is complete. 
	\noindent\subsection{ Proof of Lemma \ref{lem:bd:pot}}
	First we prove that the function $g(x):=U_{\rho,\mu}(x,\lambda) +\frac{\rho}{2}\|x-x^r\|^2_{L^{+}}$ is strongly convex with respect to variable $x$ when $2\rho\geq L_1$.  From Assumption A.1 and the  fact that $D\succeq I$ we have
		\begin{align*}
		\langle \nabla g(x_1)-\nabla g(x_2), x_1-x_2\rangle&= \langle \nabla f_\mu(x_1)-\nabla f_\mu(x_2)+\rho(A^TA+L^+)(x_1-x_2), x_1-x_2\rangle\\
		&= \langle \nabla f_\mu(x_1)-\nabla f_\mu(x_2)+2\rho D(x_1-x_2), x_1-x_2\rangle\\
		&=\langle \nabla f_\mu(x_1)-\nabla f_\mu(x_2), x_1-x_2\rangle +2\rho\|x_1-x_2\|^2_D\\
		&\geq\langle \nabla f_\mu(x_1)-\nabla f_\mu(x_2), x_1-x_2\rangle +2\rho\|x_1-x_2\|^2\\
		&\geq -L_1\|x_1-x_2\|^2+2\rho\|x_1-x_2\|^2 = (2\rho-L_1)\|x_1-x_2\|^2.
		\end{align*}
		This proves that $U_{\rho,\mu}(x,\lambda) +\frac{\rho}{2}\|x-x^r\|^2_{L^{+}}$ is strongly convex  with modulus $2\rho-L_1$.
	Using this fact, we can  bound  $U_{\rho,\mu}^{r+1}-U_{\rho,\mu}^r$ as follows:
		\begin{align}
		U_{\rho,\mu}^{r+1} - U_{\rho,\mu}^r&= U_{\rho,\mu}(x^{r+1},\lambda^{r+1}) - U_{\rho,\mu}(x^{r+1},\lambda^{r})+U_{\rho,\mu}(x^{r+1},\lambda^{r})- U_{\rho,\mu}(x^r,\lambda^r)\nonumber\\
		&\leq \frac{1}{\rho}\|\lambda^{r+1}-\lambda^r\|^2+ \langle\nabla_x U_{\rho,\mu}(x^{r+1}, \lambda^r) + \rho L^{+}(x^{r+1}-x^r),x^{r+1}-x^{r}\rangle\nonumber\\
		&-\frac{2\rho-L_1}{2}\|x^{r+1}-x^r\|^2,\label{eq:lag:bound:derive}
		\end{align}
		where the last inequity holds true due to the  strong convexity of $U_{\rho,\mu}(x,\lambda) +\frac{\rho}{2}\|x-x^r\|^2_{L^{+}}$  and \eqref{eq:due:rearg}. Now using \eqref{eq:kkt:z:app} we  have
		\begin{align}
		U_{\rho,\mu}^{r+1} - U_{\rho,\mu}^r&\le \big\langle \nabla f_\mu(x^{r+1})-G_{\mu}^{J,r},x^{r+1}-x^r\big\rangle+\frac{1}{\rho}\|\lambda^{r+1}-\lambda^r\|^2-\frac{2\rho-L_1}{2}\|x^{r+1}-x^r\|^2\nonumber\\
		&\le \frac{1}{\rho}\|\lambda^{r+1}-\lambda^r\|^2+\frac{{ L_1^2}-2\rho+L_1}{2}\|x^{r+1}-x^r\|^2+\frac{1}{2L_1^2}\|\nabla f_\mu(x^{r+1})-G_{\mu}^{J,r}\|^2\nonumber,
		\end{align}
		here the last inequality follows from \eqref{eq:rel3} with $\epsilon=L_1^2$. Taking expectation on both sides we have	
		\begin{align}
		\mathbb{E}\big[U_{\rho,\mu}^{r+1}-U_{\rho,\mu}^r\big]&\le\frac{9L_0^2(Q+4)^2}{\rho\sigma_{\min}J^2}+\frac{6L^2_1}{{\rho\sigma_{\min}}}\mathbb{E}\|x^r-x^{r-1}\|^2\nonumber\\		&+\frac{3\rho\|L^+\|}{\sigma_{\min}}\mathbb{E}\|w^r\|_{L^+}^2+\frac{L_1^2-2\rho+L_1}{2}\mathbb{E}\|x^{r+1}-x^r\|^2\nonumber\\
		&+\frac{1}{2L_1^2}\mathbb{E}\|\nabla f_\mu(x^{r+1})-\nabla f_\mu(x^{r})+\nabla f_\mu(x^{r})-G_{\mu}^{J,r}\|^2\nonumber\\
		&\le\big(\frac{9}{\rho\sigma_{\min}}+\frac{1}{L_1^2}\big)\frac{L_0^2(Q+4)^2}{J^2}+\frac{6L^2_1}{{\rho\sigma_{\min}}}\mathbb{E}\|x^r-x^{r-1}\|^2\nonumber\\
		&\hspace{-.2cm}+\frac{3\rho\|L^+\|}{\sigma_{\min}}\mathbb{E}\|w^r\|_{L^+}^2+\frac{L_1^2-2\rho+L_1+1}{2}\mathbb{E}\|x^{r+1}-x^r\|^2, \label{eq:aug:diff}
		\end{align}
		where the first inequality follows from  \eqref{eq:bd:dual} and the second inequality follows from \eqref{eq:bd_unbised}. Now we bound $V^{r+1}-V^r$.  From  the optimality condition for problem \eqref{eq:x:zero} and  the dual update \eqref{eq:lambda:zero}  we have
		\begin{align*}
		\langle &G_{\mu}^{J,r}+A^\top \lambda^{r+1} + \rho L^{+}(x^{r+1}-x^r),x^{r+1}-x\rangle\le 0,\; \forall~x\in\mathbb{R}^Q.
		\end{align*}
	Similarly, for the $(r-1)$th iteration, we have 

		\begin{align*}
		\langle &G_{\mu}^{J,r-1}+A^\top \lambda^{r}+ \rho L^{+}(x^{r}-x^{r-1}),x^{r}-x\rangle\le 0,\; \forall~x\in\mathbb{R}^Q.
		\end{align*}
		Setting $x=x^r$ in first equation, $x=x^{r+1}$ in second equation, and adding them, we obtain
		\begin{align}\label{eq:sum:kkt}
		\langle A^\top (\lambda^{r+1}-\lambda^r), x^{r+1}-x^r\rangle&\le-\langle G_{\mu}^{J,r} - G_{\mu}^{J,r-1}+ \rho L^{+}w^r, x^{r+1}-x^r\rangle.
		\end{align}
		The l.h.s can be expressed as follows:
		\begin{align}\label{eq:lhs}
		\langle  A^\top (\lambda^{r+1}-\lambda^r), x^{r+1}-x^r\rangle&=\rho\langle  Ax^{r+1}, A x^{r+1}-Ax^r\rangle\nonumber\\
		&=\frac{\rho}{2}\bigg(\|Ax^{r+1}\|^2-\|Ax^{r}\|^2+\|A(x^{r+1}-x^r)\|^2\bigg),
		\end{align}
		where the first equality follows from \eqref{eq:lambda:zero} and the second equality follows from \eqref{eq:rel1}.
	For the  r.h.s of \eqref{eq:sum:kkt} we have
		\begin{align*}
		&-\langle G_{\mu}^{J,r} - G_{\mu}^{J,r-1}+ \rho L^{+} w^r,x^{r+1}-x^r\rangle\nonumber\\
		&=-\langle G_{\mu}^{J,r} - G_{\mu}^{J,r-1}, x^{r+1}-x^r\rangle - \langle\rho L^{+} w^r, x^{r+1}-x^r\rangle \nonumber\\
		&\le  \frac{1}{2L_1}\|G_{\mu}^{J,r} - G_{\mu}^{J,r-1}\|^2+\frac{L_1}{2}\|x^{r+1}-x^r\|^2 -\rho \langle L^{+} w^r, x^{r+1}-x^r\rangle\\
		&\le \frac{3}{2L_1}\bigg(\|G^{J,r}_\mu-\nabla f_\mu(x^{r})\|^2+\|\nabla f_\mu(x^{r-1})-G^{J,r-1}_\mu\|^2+\|\nabla g_\mu(x^{r})-\nabla g_\mu(x^{r-1})\|^2\bigg)\nonumber\\
		&+\frac{L_1}{2}\|x^{r+1}-x^r\|^2 -\rho \langle L^{+} w^r, x^{r+1}-x^r\rangle,
		\end{align*}
		where the first inequality follows from \eqref{eq:rel3}. To get the second inequality we add and subtract $\nabla f_\mu(x^{r}) + \nabla f_\mu(x^{r-1})$ to  $G_{\mu}^{J,r} - G_{\mu}^{J,r-1}$ and use \eqref{eq:rel2}.   
	Taking expectation on both sides, we have
		\begin{align}\label{eq:rhs}
		&-\mathbb{E}[\langle G_{\mu}^{J,r} - G_{\mu}^{J,r-1}+ \rho L^{+}w^r, x^{r+1}-x^r\rangle]\nonumber\\
		&\le \frac{6L_0^2(Q+4)^2}{L_1J^2}+\frac{3L_1}{2}\mathbb{E}\|x^{r}-x^{r-1}\|^2+\frac{L_1}{2}\mathbb{E}\|x^{r+1}-x^r\|^2 -\rho \mathbb{E}[\langle L^{+} w^r, x^{r+1}-x^r\rangle]\nonumber\\	
		&=\frac{6L_0^2(Q+4)^2}{L_1J^2}+\frac{3L_1}{2}\mathbb{E}\|x^{r} - x^{r-1}\|^2+\frac{L_1}{2}\mathbb{E}\|x^{r+1}-x^r\|^2\nonumber\\
		&+\frac{\rho}{2}\mathbb{E}\bigg[\|x^{r}-x^{r-1}\|^2_{L^{+}}-\|x^{r+1}-x^r\|^2_{L^{+}} -\|w^r\|_{L^{+}}^2\bigg],
		\end{align}
		where the inequality follows from \eqref{eq:bd_unbised} and the last equality follows from \eqref{eq:rel1}.
	Combining \eqref{eq:lhs} and \eqref{eq:rhs}, we obtain
		\begin{align}\label{eq:rl}
		&\frac{\rho}{2}\mathbb{E}\bigg(\|Ax^{r+1}\|^2-\|Ax^{r}\|^2+\|A(x^{r+1}-x^r)\|^2\bigg)\nonumber\\
		&
		\leq \frac{6L_0^2(Q+4)^2}{L_1J^2}+\frac{3L_1}{2}\mathbb{E}\|x^{r} - x^{r-1}\|^2+\frac{L_1}{2}\mathbb{E}\|x^{r+1}-x^r\|^2\nonumber\\
		&+\frac{\rho}{2}\mathbb{E}\bigg(\|x^{r}-x^{r-1}\|^2_{L^{+}}-\|x^{r+1}-x^r\|^2_{L^{+}} -\|w^r\|_{L^{+}}^2\bigg).
		\end{align} 
		Rearranging terms in \eqref{eq:rl}, and using the definition of $V^r$ and $B$, we have
			\begin{align}\label{eq:vpart}
			\mathbb{E}[V^{r+1}-V^r]&\le \bigg(\frac{L_1}{2}+\frac{k}{2c}\bigg)\mathbb{E}\|x^{r+1}-x^r\|^2+\frac{6L_0^2(Q+4)^2}{L_1J^2}\nonumber\\
			& + \bigg(\frac{3L_1}{2}-\frac{k}{2c}\bigg)\mathbb{E}\|x^{r} - x^{r-1}\|^2-\frac{\rho}{2}\mathbb{E}\bigg(\|w^r\|_{L^{+}}^2+\|A(x^{r+1}-x^r)\|^2\bigg)\nonumber\\
			&\le \bigg(\frac{L_1}{2}+\frac{k}{2c}\bigg)\mathbb{E}\|x^{r+1}-x^r\|^2+ \bigg(\frac{3L_1}{2}-\frac{k}{2c}\bigg)\mathbb{E}\|x^{r} - x^{r-1}\|^2\nonumber\\
			& \qquad-\frac{\rho}{2}\mathbb{E}\|w^r\|_{L^{+}}^2 +\frac{6L_0^2(Q+4)^2}{L_1J^2}.
			\end{align}
			Recall the definition of $P^{r+1}:=U_{\rho,\mu}^{r+1}+ cV^{r+1}$. 
		Utilizing \eqref{eq:aug:diff} and \eqref{eq:vpart} eventually we obtain	
			\begin{align*}
			\mathbb{E}\big[P^{r+1}-P^{r}\big]&\le -\alpha\mathbb{E}\|x^{r+1} - x^r\|^2 -\alpha_2\mathbb{E}\|w^r\|_{L^{+}}^2 + \alpha_3\frac{L_0^2(Q+4)^2}{J^2},
			\end{align*}
			where we have  
				\begin{align*}
				&\alpha_1:= \rho^2-(2cL_1+L_1/2+L^2/2+1/2)\rho-\frac{6L_1^2}{\sigma_{\min}},\nonumber\\
				& \alpha_2:=\frac{3\rho\|L^{+}\|}{\sigma_{\min}}-\frac{c \rho}{2},\; \alpha_3=\frac{9}{\rho\sigma_{\min}}+\frac{6c+1}{L_1}.
				\vspace{-.3cm}
				\end{align*}
			The lemma is proved.

		\noindent\subsection{ Proof of Lemma \ref{lemma:lower:bound}}
		\vspace{-.1cm}
		From \eqref{eq:kkt:z:app} we have
		\begin{align}
		\|A^\top\lambda^{r+1}\|^2\leq 2\|G_{\mu}^{J,r}\|^2+2\rho^2 \|L^{+}(x^{r+1}-x^r)\|^2.
		\end{align}
		Recall that $\sigma_{\min}$ is the smallest nonzero eigenvalue of $A^TA$. Also from the fact that $\lambda^0=0$
		we have that the dual variable  lies in the column space of $A$, thus we get that
		\begin{align}\label{lam}
		\|\lambda^{r+1}\|^2\leq \frac{2}{\sigma_{\min}}\|G_{\mu}^{J,r}\|^2+\frac{2\rho^2}{\sigma_{\min}}\|L^+(x^{r+1}-x^r)\|^2.
		\end{align}
		From the definition of the potential function we have
		{\small
			\begin{align}\label{pot}
			P^{r+1} &= f_\mu(x^{r+1}) +\frac{\rho}{2}\|Ax^{r+1}+\frac{1}{\rho}\lambda^{r+1}\|^2-\frac{1}{2\rho}\|\lambda^{r+1}\|^2\nonumber\\
			&+\frac{c\rho}{2}\|Ax^{r+1}\|^2 +\frac{c\rho}{2}\|x^{r+1}-x^r\|_B^2,
			\end{align}}where $B:=L^++\frac{k}{c\rho}I$. From Assumption A.2 we have that $f$ is lower bounded, therefore,  $f_\mu$ is lower bounded too, i.e. there exists $\underline{f}$ such that $\underline{f}\leq f_\mu(x)$ for all $x\in\dom(f_\mu)$. Plugging \eqref{lam} in \eqref{pot}, and utilizing the fact that  $\|Ax^{r+1}+\frac{1}{\rho}\lambda^{r+1}\|^2\geq 0$  we obtain
			\begin{align*}
			P^{r+1} &\geq \frac{-1}{\rho\sigma_{\min}}\|G_{\mu}^{J,r}\|^2-\frac{\rho}{\sigma_{\min}}\|L^{+}(x^{r+1}-x^r)\|^2+\frac{c\rho}{2}\|x^{r+1}-x^r\|_B^2+\underline{f}\\
			&\geq \frac{-1}{\rho\sigma_{\min}}\|G_{\mu}^{J,r}\|^2 + \frac{\rho}{\sigma_{\min}}\|x^{r+1}-x^r\|^2_{\frac{c\sigma_{\min}}{2}L^+-(L^+)^2}+ \underline{f},
			\end{align*}
			where the last inequality is due to the fact that $\frac{k}{2}\|x^{r+1}-x^r\|^2\geq0$. Notice that $L^+$ is a symmetric PSD matrix. Therefore, picking constant $c$ large enough such that 
			$c\geq \frac{2\|L^+\|}{\sigma_{\min}}$, 
			we have $\frac{c\sigma_{\min}}{2}L^+-(L^+)^2\succeq 0$. Hence, with this choice of $c$ we get the following bound for the potential function
		\begin{align}
		P^{r+1} \geq -\frac{1}{\rho\sigma_{\min}}\|G_{\mu}^{J,r}\|^2 + \underline{f}.
		\end{align}
		Taking expectation on both sides we have
		\begin{align}\label{pot2}
		\mathbb{E}[P^{r+1}] &\geq -\frac{1}{\rho\sigma_{\min}}\mathbb{E}\|G_{\mu}^{J,r}\|^2 + \underline{f}\nonumber\\
		&\geq -\frac{L_0(Q+4)^2}{\rho\sigma_{\min}J^2} + \underline{f},
		\end{align}
		where the last inequality follows from \eqref{eq:nes_them4}. To complete the proof we only need  to set $\underline{P}=-\frac{L_0(Q+4)^2}{\sigma_{\min}J^2} + \underline{f}$.
	
	\noindent\subsection{ Proof of Theorem \ref{thm:conv}}
	First we bound the stationarity gap given in \eqref{eq:opt_gap}. We have
		\begin{align}
		\|\nabla_x U_{\rho,\mu}(x^{r+1}, \lambda^r)\|^2 &= \|\nabla f_\mu(x^{r+1})+A^\top \lambda^{r+1}\|^2\nonumber\\
		&=\|\nabla f_\mu(x^{r+1})-G_{\mu}^{J,r}-\rho L^{+}(x^{r+1}-x^r)\|^2\nonumber\\
		&\le 2\|\nabla f_\mu(x^{r+1})-G_{\mu}^{J,r}\|^2+2\rho^2\|L^{+}(x^{r+1}-x^r)\|^2\nonumber\\
		&\leq 2\|\nabla f_\mu(x^{r+1})-\nabla f_\mu(x^{r})+\nabla f_\mu(x^{r})-G_{\mu}^{J,r}\|^2+2\rho^2\|L^{+}(x^{r+1}-x^r)\|^2\nonumber\\
		&\leq 4\|\nabla f_\mu(x^{r+1})-\nabla f_\mu(x^{r})\|+4\|\nabla f_\mu(x^{r})-G_{\mu}^{J,r}\|^2\nonumber\\
		&+2\rho^2\|L^{+}(x^{r+1}-x^r)\|^2,\nonumber
		\end{align}
		where the first equality follows from \eqref{eq:lambda:zero}, the second equality follows from \eqref{eq:kkt:z:app}, and the first inequality follows from \eqref{eq:rel2}. Taking expectation on both sides gives
		\begin{align}
		\mathbb{E}\|\nabla_x U_{\rho,\mu}(x^{r+1}, \lambda^r)\|^2	&\leq 4L_1^2\mathbb{E}\|x^{r+1}-x^{r}\|+\frac{4L_0^2(Q+4)^2}{J^2}+2\rho^2\mathbb{E}\|L^{+}(x^{r+1}-x^r)\|^2, \label{eq:bd:grad2}
		\end{align}
		where the inequality follows from \eqref{eq:bd_unbised}. Next, we bound the expected value of the constraint violation. Utilizing the equation  \eqref{eq:lambda:zero} we have
		\begin{align*}
		\|Ax^{r+1}\|^2=\frac{1}{\rho^2}\|\lambda^{r+1}-\lambda^r\|^2.
		\end{align*}
		Taking expectation and utilizing \eqref{eq:mu:difference:bound} we reach
		\begin{align}\label{eq:bd:vio}
		\mathbb{E}\|Ax^{r+1}\|^2&=\frac{1}{\rho^2}\mathbb{E}\|\lambda^{r+1}-\lambda^r\|^2\leq \frac{9L_0^2(Q+4)^2}{\rho^2\sigma_{\min}J^2}\nonumber\\
		&+\frac{6L^2_1}{\rho^2\sigma_{\min}}\mathbb{E}	\|x^r-x^{r-1}\|^2+\frac{3\|L^{+}\|}{\sigma_{\min}}\mathbb{E}\|w^r\|_{L^+}^2.	\end{align}
		Summing  \eqref{eq:bd:grad2} and \eqref{eq:bd:vio}, we have the following bound for the stationarity gap 
		\begin{align}
		\Phi_\mu^{r+1}&\le\beta_1\mathbb{E}\|x^{r+1}-x^r\|^2+\beta_2\mathbb{E}\|x^r-x^{r-1}\|^2+\beta_3\mathbb{E}\|w^r\|_{L^{+}}^2+\beta_4\frac{L_0^2(Q+4)^2}{J^2},\label{eq:Q:bd}
		\end{align}
		where $\beta_1, \beta_2, \beta_3, \beta_4$ are positive constants given by 
		\begin{align*}
		&\beta_1=4{L_1}^2+2\rho^2{\|L^{+}\|^2}, \; \beta_2=\frac{6{L_1}^2}{\rho^2\sigma_{min}},\\
		&\beta_3=\frac{3\|L^{+}\|}{\sigma_{min}},\; \beta_4 = \frac{9+4\rho^2\sigma_{\min}}{\rho\sigma_{\min}}.
		\end{align*} 
		Summing both sides of \eqref{eq:Q:bd} over $T$ iterations,  we get
		\begin{align}
		\sum_{r=1}^{T}\Phi_\mu^{r+1}&\le \sum_{r=1}^{T-1}(\beta_1+\beta_2)\mathbb{E}\|x^{r+1}-x^{r}\|^2+\sum_{r=1}^{T} \beta_3\mathbb{E}\| w^r\|_{L^{+}}^2\nonumber\\
		&+\beta_2\mathbb{E}\|x^1-x^0\|^2 +\beta_1\mathbb{E}\|x^{{T}+1}-x^{T}\|^2+T\beta_4\frac{L_0^2(Q+4)^2}{J^2}.\label{eq:sum:Q}
		\end{align}
		Applying Lemma \ref{lem:bd:pot} and  summing both sides of \eqref{eq:diff:bd}  over $T$ iterations, we obtain
		\begin{align}
		\mathbb{E}\big[P^1-P^{T+1}\big]&\ge \sum_{r=1}^{T-1}{\alpha_1}\mathbb{E}\|x^{r+1}-x^{r}\|^2+\sum_{r=1}^{T} \alpha_2\mathbb{E}\| w^r\|_{L^+}^2\nonumber\\
		&+\alpha_1\mathbb{E}\|x^{T+1}-x^T\|^2-T\alpha_3\frac{L_0^2(Q+4)^2}{J^2}.\label{eq:sum:P}
		\end{align}
		Let us set $\tau=\frac{\max(\beta_1+\beta_2, \beta_3)}{\min(\alpha_1, \alpha_2)}$. Combining the two inequalities \eqref{eq:sum:Q} and \eqref{eq:sum:P} and utilizing the fact that  $\mathbb{E}[P^{T+1}]$ is lower bounded by $\underline{P}$, we arrive at the following inequality
		\begin{align}\label{eq:sum:Q:2}
		\sum_{r=1}^{T}\Phi_\mu^{r+1}&\le  \tau\mathbb{E}[P^1-\underline{P}]+\beta_2\mathbb{E}\|x^1-x^0\|^2+T\big(\tau \alpha_3+\beta_4\big)\frac{L_0^2(Q+4)^2}{J^2}.
		\end{align}
		Because $u$ is a uniform random number from $\{0,1,\cdots, T-1\}$ we have
		\begin{align}\label{eq:exp:J}
		\mathbb{E}_u[\Phi_\mu^u]=\frac{1}{T}\sum_{r=1}^{T}\Phi_\mu^{r+1}.
		\end{align}
		Combining \eqref{eq:sum:Q:2} and \eqref{eq:exp:J} implies the following
		\begin{align*}
		\mathbb{E}_u[\Phi_\mu^u]&\le \frac{\tau\mathbb{E}[P^1-\underline{P}]+\beta_2\mathbb{E}\|x^1-x^0\|^2}{T}+ \big(\tau \alpha_3+\beta_4\big)\frac{L_0^2(Q+4)^2}{J^2}.
		\end{align*}
		Setting 
		\begin{align}
		&\gamma_1=\tau\mathbb{E}[P^1-\underline{P}]+\beta_2\mathbb{E}\|x^1-x^0\|^2,\;\gamma_2=\tau \alpha_3+\beta_4,
		\end{align}
		we conclude the proof.
	
		\bibliographystyle{IEEEbib}
		\bibliography{ref_davood}

\end{document}